
%
%
%
%
%
%
%
%
\magnification =1200
\footline={\hss\tenrm\folio\hss}
\normalbaselineskip = 1.2\normalbaselineskip
\normalbaselines

\font\titlefont=cmbx10 scaled \magstep4
\font\authorfont=cmbx10 scaled \magstep2

\def\hat{\widehat}

\def\Z{{\bf Z}}

\def\Rat{{Rat}}
\def\Ind{{Ind}}
\def\R{{\bf R}}

\def\C{{\bf C}}

\def\Spin{{\rm Spin}}

\def\and{\hbox{\quad and \quad}}

%
%

\nopagenumbers

\null
\vskip 1.5 cm
\centerline{\titlefont A note on the index bundle over}
\vskip .5 cm
\centerline{\titlefont  the moduli space of monopoles}
\vskip 2 cm

\centerline{\authorfont John D. S. Jones$^1$}

\vskip 1.5 cm
\centerline{\authorfont Michael K. Murray$^2$}

\vskip 1.5 cm

\vfill

\centerline{\bf Preprint: 15 July 1994}

\vfill

\noindent 1. Mathematics Institute, University of Warwick, Coventry   CV 4
7AL, United Kingdom.
{\it jdsj@maths.warwick.ac.uk}

\noindent 2. Department of Pure Mathematics, The University of Adelaide,
Adelaide, SA 5005 Australia. {\it mmurray@maths.adelaide.edu.au}

\vfill

\eject
\nopagenumbers

\vskip 1cm

\noindent{\bf Abstract.}  Donaldson has shown that the moduli
space of monopoles $M_k$ is diffeomorphic to the space   $\Rat_k$ of based
rational maps  from the two-sphere to itself.
We use this diffeomorphism to give an explicit description
of the bundle on $\Rat_k$ obtained by pushing out
the index bundle from $M_k$. This gives an alternative and more explicit
proof of some earlier results of Cohen and Jones.


\vfill\eject

\pageno=1
\footline={\hss\tenrm\folio\hss}

\noindent{\bf  1. Introduction \hfill}

In [4] Cohen and Jones study the topological type of the index bundle
of various families of Dirac operators arising in the theory of
monopoles and the relation between these index bundles and
representations of the braid groups. The methods used in this general
study were those of algebraic topology and index theory. For example,
it was shown that using Donaldson's diffeomorphism  between monopoles
and based rational maps [5] and the relation between the space of
based rational maps and the braid group that the $K$-theory class of
the index bundle over the space of monopoles is completely determined
by representations of the braid groups. In this note we show how
Donaldon's diffeomorphism gives rise to a  simple explicit
characterisation of the bundle over the space of rational maps
corresponding to this index bundle. The corresponding representation
of the braid group is readily identified.

In more detail, let $M_k$ denote the moduli space of framed monopoles
of charge $k$ over $\R^3$ with  structure group $SU(2)$. Donaldson in
[5] defines an explicit  diffeomorphism
$$
M_k \to \Rat_k.
\eqno(1.1)
$$
from $M_k$ to the space of all based rational maps of degree $k$ from
the two-sphere to itself. It can be shown [10] that if $(A, \Phi)$ is
a monopole then the space of $L^2$ solutions of the Dirac equation
coupled to $(A, \Phi)$ has dimension $k$ where $k$ is the charge of
the monopole. This defines a complex vector bundle  over $M_k$ which,
in fact, has a hermitian inner product and a real structure and hence
has structure group $O(k)$, the group of real, orthogonal, $k$ by $k$
matrices. We denote by $\Ind_k$,  the corresponding principal $O(k)$
bundle on $M_k$. This is the index bundle and the point of this note
is to show that we can explicitly    describe  the  bundle over
$\Rat_k$ which is its `push-out' under the diffeomorphism (1.1). In
fact this bundle arises quite naturally in Donaldson's work. We show
further, that over the space $\Rat^0_k$ of rational maps with distinct
poles, this bundle has a reduction to the finite subgroup of $O(k)$ of
signed permutations.  This bundle therefore corresponds to a
homomorphism of $\pi_1(\Rat^0_k)$ into the group of signed
permutations. This homomophism can be readily calculated and it is,
not suprisingly, the same as that in [4].

\bigskip

\noindent{\bf  2. Monopoles and the index bundle \hfill}

The index bundle $\Ind_k$ is defined over $M_k$ the moduli space of
framed monopoles of charge $k$.  To define this moduli space  consider
first  pairs $(A, \Phi)$ consisting of an $SU(2)$ connection $A$ on
$\R^3$ and an $su(2)$ `Higgs field' $\Phi\colon \R^3 \to su(2)$, where
$su(2)$ is the Lie algebra of $SU(2)$.  The Bogomolny equations for
such a  pair are
$$
\star F_A = \nabla_A\Phi
\eqno(2.1)
$$
where $F_A$ is the curvature of $A$, $\nabla_A\Phi$ the covariant
derivative of $\Phi$ with respect to $A$ and $\star$ is the Hodge dual
on forms. To be a monopole $(A, \Phi)$ has to satisfy the Bogomolny
equations and also  certain boundary conditions that we do not need to
detail here. We refer to the book [2] as a good general reference for
these and other details.  However we do need to note that one of the
boundary conditions is that the Higgs field gives rise to a map
$\Phi^\infty$  from the two-sphere at infinity in $\R^3$ and that this
takes values on the two-sphere inside $su(2)$. The degree of this map
is called the magnetic charge of the monopole and we shall denote it
by $k$. To be a framed monopole we require that
$$
\lim_{t \to 0} \Phi(0, 0, t) = \pmatrix{ i & 0 \cr 0 & -i \cr}.
\eqno(2.2)
$$
The group of gauge transformations, that is  smooth maps $g \colon \R^3
\to SU(2) $, preserving the boundary conditions, acts on  pairs $(A,
\Phi)$ to give  new pairs $(A^g, \Phi^g)$ defined by
$$
A^g = g^{-1}Ag + g^{-1}dg  \quad \rm{and} \quad \Phi^g = g^{-1} \Phi g.
\eqno(2.3)
$$
These gauge transformations preserve the set of solutions of the
Bogomolny equations (2.1).  They also  preserve the framing condition
(2.2) if  the limit of $g(0,0,t)$ as $t$ goes to infinity is diagonal.
A framed gauge transformation is defined to be one such that
$$
\lim_{t \to 0} g(0, 0, t) = 1.
\eqno(2.4)
$$
The group of framed gauge transformations acts freely on the set of
framed monopoles and the quotient is $M_k$ the moduli space of framed
monopoles. It is a manifold of dimension $4k$ [11].

Given a pair $(A,\Phi)$ we can form the coupled Dirac operator $D_0$
acting on sections of the trivial spinor bundle over $\R^3$ with fibre
$\C^2 \otimes \C^2$:
$$
D_0 = \sum_{i= 1}^3 \sigma_i\otimes {\nabla_A}_i - 1\otimes \Phi
\eqno(2.5)
$$
where the $\sigma_i$ are the matrices defining the action of the
Clifford algebra of $\R^3$ on $\C^2$. It is known [10] that
the space of $L^2$-solutions of the equation $D_0 \psi = 0$,
satisfying the given boundary conditions, has dimension $k$ where
$k$ is the charge of the monopole. The group of framed gauge
transformations acts on the  spinor bundle and on
the Dirac operator by conjugation  hence quotienting gives rise to
a vector bundle over $M_k$.

The Dirac operator $D_0$  acts on sections of a trivial bundle with
fibre $\C^2 \otimes \C^2$ where one factor is acted on by the group
$\Spin(3) = SU(2)$ and the other is acted on by the $SU(2)$ from the
monopole bundle. In both cases the structure group is $SU(2) = Sp(1)$
and hence the individual bundles have quaternionic structures.
Therefore there is a real structure on the tensor product, that is a
conjugate linear map $r$ such that  $r^2 = 1$. The space $L^2(\R^3,
\C^2\otimes \C^2)$ has a hermitian inner product defined by integrating
the hermitian inner product on $\C^2\otimes \C^2$ and this restricts
to a hermitian inner product on the kernel of $D_0$. The real
structure  map $r$ preserves this natural hermitian  inner product
$\langle\ ,\ \rangle$ and therefore defines an orthogonal form by $(v,
w) = \langle v, r(w)\rangle$. It makes sense therefore to consider the
space of all real, orthonormal frames for the kernel of $D_0$. Here
real means fixed by the real structure. This space is acted on freely
by $O(k)$ the group of real, orthogonal, $k$ by $k$ matrices. Hence we
have constructed a principal $O(k)$ bundle over $M_k$ which we shall
call the index bundle and denote by $\Ind_k$.

\bigskip
\noindent{\bf  3. The ADHMN construction \hfill}

The ADHM construction for instantons [1] as generalised by Nahm to
monopoles [10] associates to every pair $(A, \Phi)$ satisfying the
Bogomolny equations and the appropriate monopole bondary conditions  a
rank $k$ bundle $N$ over the interval $(-1, 1) \subset \R$ equipped
with a connection $\nabla$  and three bundle endomorphisms $T_i$. If
we trivialise the bundle with covariantly constant sections then the
$T_i$ become matrices satisfying Nahm's equations:
$$
\matrix{
{dT_1 \over dz} = [T_2, T_3]\cr
{dT_2 \over dz} = [T_3, T_1]\cr
{dT_3 \over dz} = [T_1, T_2]\cr
}
\eqno(3.1)
$$
and some boundary conditions.  Again it is not important precisely what
the boundary conditions are and we refer the reader to [2] for
details. We do need to note however that the $T_i$ are analytic and
have simple poles at $\pm 1$. Let $t_i$ denote the residues of the
$T_i$ at $-1$. It follows from Nahm's equations (3.1) that the
residues are a representation of $su(2)$. It is one of the boundary
conditions that they must  in fact  be irreducible, and hence  $-it_3$
has eigenvalues $-(k-1), \dots, (k-1)$. Although it was not explicit
in Nahm's work one can follow through the constructions in [7] to see
that the framing of the monopole means that we also have given a unit
vector $v$ in the $(k-1)$ eigenspace of $-iT_3$.

The connection with the index bundle follows from the fact that
fibre of the bundle $N$ over the point $t \in (-1, 1)$ is the $L^2$
kernel of the coupled Dirac operator
$$
D_t =  \sum_{i= 1}^3 \sigma_i\otimes  (\nabla_A)_i - 1\otimes (\Phi + it ).
\eqno(3.2)
$$
The hermitian and real structures of the $\C^2\otimes \C^2$ bundle over
$\R^3$ pass to the bundle $N$ as follows. Firstly the vector space
$N_t$ is a subspace of the hermitian inner product space $L^2(\R^3,
\C^2\otimes \C^2)$ so it inherits a hermitian inner product by
restriction. Secondly the  real structure is conjugate linear so it
maps an element of the kernel of $D_t$ to an element in the kernel of
$D_{-t}$. Hence it defines a conjugate linear map from $N_t$ to
$N_{-t}$. There is an orthogonal projection
$$
\pi_t \colon L^2(\R^3, \C^2\otimes \C^2) \to N_t
\eqno(3.3)
$$
and this is used to define the connection and the endomorphisms $ T_i$
by
$$
\nabla(\chi) = \pi({d\chi\over dt}) \quad {\rm and}\quad
T_i(\chi) = \pi(x^i\chi).
\eqno(3.4)
$$

We can identify  $N_0$ with $\C^k$ by choosing an orthonormal frame
consisting of real vectors and using parallel transport to  extend
this to a frame at every point of $(-1, 1)$ and hence make the
$T_i$ into matrices. So to every triple $(A,  \Phi, \{\psi_a\})$,
where  $(A, \Phi)$ is a monopole and $\{\psi_a\}$ is an orthonormal
basis of real vectors in the space $N_0$ of $L^2$ solutions of the
Dirac equation coupled to $(\nabla,\Phi)$, we have associated three
matrix valued functions $T_i$ on the interval $(-1, 1)$ satisfying
Nahm's equations (3.1) and a $v$ as above.   From the definition (3.4)
we see that  if  we gauge transform the triple $(A,  \Phi,
\{\psi_a\})$ then the $T_i$  are left unchanged. The same is also true
of $v$ because we are using framed gauge transformations. If the basis
$\{\psi_a\}$ is multiplied by an element of $O(k)$ then the $T_i$ are
conjugated by this same element and the $v$ is multiplied by it.  So
the  image of the index bundle after pushing out with the ADHMN
construction is all  $(T_1, T_2, T_3, v)$ where the $T_i$ are a
solution to Nahm's equations with appropriate boundary conditions and
the $v$ is defined as above. The action of $O(k)$ on this principal
bundle is conjugation of the $T_i$ and multiplication of the $v$.

\bigskip \noindent{\bf 4. The index bundle over $\Rat_k$ \hfill}

Donaldson shows that the space $N_k$, and hence $M_k$, is
diffeomorphic to the space $\Rat_k$ of based rational maps from the
two sphere to the two sphere of degree $k$. The diffeomorphism  is
defined as follows.  Let us first combine Nahm's matrices as
$$
A_0 = T_1 + iT_2 ,\quad A_1 = -i T_3, \and A_2 = T_1 - iT_2.
\eqno(4.1)
$$
Then consider solutions $u \colon (-1, 1) \to \C^k$ of the
differential equation
$$
{d u \over ds } - {1 \over 2} A_1 u = 0.
\eqno(4.2)
$$
There is a unique solution $u$ with the property that
$$
\lim_{s \to -1} s^{-(k-1)/2} u(s) = v.
\eqno(4.3)
$$
Define $B = -A_0(1)$ and $W = u(1)$.  Then Donaldson shows that $B$ is
a symmetric matrix and that $W$ is a cyclic vector for $B$; that is
$\{W, BW, \dots, B^{k-1}W\}$ are linearly independent. He also shows
that the space $\hat{\Rat}_k$ of all such pairs $(B, W)$ with $B$
symmetric and $W$ cyclic for $B$ is a principal $O(k)$ bundle over
$\Rat_k$. The projection $\hat{\Rat}_k \to \Rat_k$ is given by
$$
(B, W) \mapsto f(z) = W^t(zI - B)^{-1} W,
\eqno(4.4)
$$
and the $O(k)$ action is given by conjugation on $B$ and left
multiplication on $W$. The map $(T_1, T_2, T_3, p) \mapsto (B, W)$ is
equivariant with respect to the $O(k)$ actions and Donaldson's result
shows that it descends to a diffeomorphism from $N_k$ to $\Rat_k$. It
follows that the principal $O(k)$ bundle $\hat{\Rat}_k \to \Rat_k$ is
the pushout of the index bundle under Donaldson's diffeomorphism from
$M_k$ to $\Rat_k$.

\bigskip \noindent{\bf  5. The reduction over $\Rat^0_k$ \hfill}

The space $\Rat^0_k$ is the space of based rational maps with distinct
poles. Consider  a diagonal matrix $B$ with distinct entries $(b_1,
\dots, b_k)$ and a vector $W = (W_1, \dots, W_k)$ with non-zero
components. It follows from the Vandermonde determinant that $W$ is
cyclic for $B$. Denote by $\hat{\Rat}_k^0$ the set of all such $(B,
W)$. The rational map defined by the pair $(B,W)$ using the projection
$\hat{\Rat}_k\to \Rat_k$, is
$$
f(z) = \sum_{i=1}^k { W_i^2 \over z - b_i}
\eqno(5.1)
$$
so that $\hat{\Rat}_0$ covers the space $\Rat_k^0$ of all
rational maps with distinct poles.  Moreover, this set is stable
under the action of the group $\Sigma_k^{\pm}$ of all signed
permutations matrices, that is the subgroup of $O(k)$ generated by
the diagonal matrices with  plus or minus one on the diagonal and
the permutation matrices.  Indeed $\hat{\Rat}_k^0$ is a principal
$\Sigma_k^{\pm}$ bundle over $\Rat_k^0$ and it defines a reduction
of the restriction of the bundle bundle $\hat{\Rat}_k$ to
$\Rat_k^0$ to a $\Sigma_k^{\pm}$ bundle.  It follows that over
$\Rat_k^0$ the bundle $\Ind_k$ has a reduction to $\Sigma_k^{\pm}$
and therefore it defines a homomorphism of $\pi_1(\Rat_k^0)$ to
$\Sigma_k^{\pm}$.
In [3] it is shown that $\pi_1(\Rat_k^0)$ is the
semi-direct product of the braid group on $k$ strings $\beta_k$
and $\Z^k$ where $\beta_k$ acts on $\Z^k$ via the natural
homomorphism to $\Sigma_k$.  This group maps naturally onto
$\Sigma_k^{\pm}$ which is the semi-direct product of the symmetric
group $\Sigma_k$ and $(\Z/2)^k$.  By considering generators of
$\pi_1(\Rat_k^0)$ it is possible to show that the reduction of the
bundle above corresponds exactly to this homomorphism; for the
details see [3].  The topological implications
of this fact are given in [4].

\vfill\eject
\noindent{\bf  6. The Dirac equation \hfill}

It is interesting to consider what this reduction means for the
solutions of the Dirac equation. The reduction of the principal
bundle $\Ind_k$ to $\Sigma_k^{\pm}$ corresponds to extra structure
on the vector bundle  of solutions of the Dirac equation.
To see what this is note that an orthonormal basis of solutions
of the Dirac equation would correspond to a reduction to the
identity subgroup. This is more than we actually have.  However
if we choose not a basis of vectors but an (unordered) set of $k$
orthogonal lines that
span the space then it is easy to see that this gives rise to
an orthonormal basis up to the action of $\Sigma_k^{\pm}$.

So the reduction we have constructed corresponds to  being able to find
a  set of $k$ orthogonal lines in the space of all solutions of the
Dirac equation.  One way of doing this has already been noted in [9].
This uses the fact that for widely seperated monopoles there are
solutions of the Dirac equation concentrated about each of the
monopole locations. The reduction that we have given appears to be
different to this. It can be understood by reference to Hitchin's
twistor construction of solutions of Nahm's equations [7]. We shall
sketch here, without proof how this occurs,  and refer the reader to
[7] and [8] for details. Recall that Hitchin showed in  [6] that a
monopole is determined by a certain algebraic curve $S$ in $TS^2$ the
tangent bundle of the two sphere. This is the so-called spectral
curve. The role played by $TS^2$ is that it parametrises the set of
all oriented lines (not necessarily through the origin) in $\R^3$. The
points of the two sphere correspond to the direction of the line and
the fibres of the projection $TS^2 \to S^2$ correspond to the families
of parallel lines. Denote by $F$ the fibre of all lines in the $x^3$
direction. In  subsequent  work [7] Hitchin  showed  that the space of
solutions of the Dirac equation, $N_z$  can be identified with the
space
$$
H^0(S, L(k-1))
\eqno(6.1)
$$
of holomorphic sections of a certain line bundle over $S$. Which line
bundle is not important for this discussion  and we refer to Hitchin's
papers [6] and [7] for the details. Hitchin also proves that the
restriction map of a section to the intersection of the fibre $F$ with
the curve $S$ is an isomorphism. Generically this intersection is $k$
distinct points. In fact these $k$ points are the poles of the
rational map [8]. In such a case we can define a line in $N_z$ by
considering those sections that vanish on restriction to all but one
of the points. By changing the point we generate $k$ lines and these
span the space $N_z$. It was shown by  Hurtubise in [8] that  these
lines are orthogonal and that if they are chosen as an orthogonal
basis  then they determine a $B$ that is diagonal and a $W$ with
non-zero components. It follows that the reduction we have described
in terms of the space of rational maps corresponds to sections whose
restriction to the fibre $F$ are supported at just one point. It would
be interesting to understand what this means for solutions of the
Dirac equation in $\R^3$. This would mean unravelling the twistor
correspondence in more detail.

\bigskip
\noindent{\bf  Acknowledgements \hfill}

MKM thanks the Mathematics Research Centre, University of
Warwick for hospitality and support.

\bigskip
\noindent{\bf  References \hfill}

\baselineskip=15pt
\parskip = 0 pt

\item{[1]}
Atiyah, M.F., Drinfeld, V.G., Hitchin, N.J. and Manin, Yu.I.:
{\it Construction of instantons},
Phys. Letts. {\bf 65A},
185--187,
1978.

\item{[2]}
Atiyah, M.F. and Hitchin, N.J.:
{\it The Geometry and Dynamics of Magnetic  Monopoles},
Princeton University Press, Princeton,
1988.

\item{[3]}
Cohen, R.L. and Jones, J.D.S.:
{\it Representations of braid groups and
operators coupled to monopoles},
in {\it Geometry of low-dimensional manifolds: 1},
Donaldson S.K. and Thomas C.B. editors,
London Mathematical Society Lecture Notes {\bf 150},
191--205,
1990.

\item{[4]}
Cohen, R.L. and Jones, J.D.S.:
{\it Representations of the braid group and
operators coupled to monopoles},
Commun. Math. Phys. {\bf 158},
241--266,
1993.

\item{[5]}
Donaldson, S.K.:
{\it Nahm's equations and the classification of
monopoles},
Commun. Math. Phys. {\bf 96},
387--407,
1984.

\item{[6]}
Hitchin, N.J.:
{\sl Monopoles and geodesics.}
Commun. Math. Phys. {\bf 83},
579--602
(1982).

\item{[7]}
Hitchin, N.J.:
{\sl On the construction of monopoles. }
Commun. Math. Phys. {\bf 89},
145--190,
(1983).

\item{[8]}
Hurtubise, J.C.:
{\sl Monopoles and rational maps: A note on a theorem of Donaldson.}
Commun. Math. Phys. {\bf 100},
191--196
(1985).

\item{[9]}
Manton, N.S. and Schroers, B.J.:
{\sl Bundles over moduli spaces and the quantisation of BPS
monopoles},
Annals of Physics,
{\bf 225},
290--338,
(1993).

\item{[10]}
Nahm, W.
{\it The construction of all self-dual
monopoles by the ADHM method},
in {\it Monopoles in Quantum Field Theory},
Craigie, N., et al editors,
World Scientific, Singapore,
1982.

\item{[11]}
Taubes, C.H.:
{\sl Monopoles and Maps from $S^2$ to $S^2$;
                        the topology of the Configuration Space.}
Commun. Math. Phys. {\bf 95},
345--391
(1984).

\bye